\documentclass[a4paper,12pt,reqno]{amsart}
\usepackage{latexsym}
\usepackage{amsmath,amsthm,amssymb,amscd}
\usepackage{fullpage}
\usepackage{xcolor}
\usepackage{parskip}
\usepackage{mathtools}
\mathtoolsset{showonlyrefs}

\usepackage[activate={true,nocompatibility},final,tracking=true,kerning=true,spacing=true]{microtype}
\usepackage[pdfencoding=unicode,pdftex,bookmarks=true,bookmarksnumbered]{hyperref}
\definecolor{citecolour}{rgb}{0.0, 0.0, 0.8}
\colorlet{linkcolour}{green!50!black}
\hypersetup{colorlinks,breaklinks,
		linkcolor=linkcolour,citecolor=citecolour,
		filecolor=linkcolour, urlcolor=linkcolour}
\usepackage{version}

\newtheorem{theorem}{Theorem}

\newtheorem*{proposition*}{Proposition A}
\newtheorem*{proposition**}{Proposition B}
\newtheorem*{proposition***}{Proposition C}
\newtheorem*{proposition****}{Proposition D}
\newtheorem*{proposition*****}{Proposition G}
\newtheorem*{proposition******}{Proposition H}

\newtheorem{lemma}[theorem]{Lemma}
\newtheorem*{corollary*}{Corollary D}
\newtheorem*{corollary**}{Corollary E}
\newtheorem*{corollary***}{Corollary F}

\newtheorem{corollary}[theorem]{Corollary}

\newtheorem*{lemma*}{Lemma}
\newtheorem*{theoremA}{Theorem A}
\newtheorem*{theoremC}{Theorem C}
\theoremstyle{remark}

\newtheorem{remark}[theorem]{Remark}

\numberwithin{equation}{section}
\usepackage[shortlabels]{enumitem}
\begingroup
    \makeatletter
    \@for\theoremstyle:=definition,remark,plain\do{%
        \expandafter\g@addto@macro\csname th@\theoremstyle\endcsname{%
            \addtolength\thm@preskip\parskip
            }%
        }

\DeclareMathOperator{\F}{\mathfrak{F}}
\DeclareMathOperator{\A}{\mathfrak{A}}
\DeclareMathOperator{\X}{\mathfrak{X}}

\DeclareMathOperator{\Zent}{\mathbf{Z}}
\DeclareMathOperator{\Soc}{\mathbf{Soc}}
\DeclareMathOperator{\Oh}{\mathbf{O}}

\DeclareMathOperator{\cent}{\mathbf{C}}

\DeclareMathOperator{\fitt}{\mathbf{F}}
\DeclareMathOperator{\fratt}{\mathbf{\Phi}}
\DeclareMathOperator{\Aut}{\mathbf{Aut}}
\DeclareMathOperator{\fast}{\mathbf{F}^{\ast}}
\DeclareMathOperator{\Lay}{\mathbf{E}}

\newcommand{\R}{\mathbb{\R}}

\includeversion{comment}
\renewcommand{\leq}{\leqslant}
\renewcommand{\geq}{\geqslant}

\newenvironment{proofof}{{\bf {Proof.} }}{\hfill $\blacksquare$ \\}

\begin{document}
\title{Large subgroups in finite groups}
\author[S. Aivazidis and 
T.\,W. M\"uller]{Stefanos Aivazidis$^{\dagger}$ and
Thomas M\"uller$^*$} 

\address{$^{\dagger}$Einstein Institute of Mathematics, Edmond J. Safra Campus, The Hebrew University of Jerusalem, Givat Ram, Jerusalem 9190401, Israel.}

\address{$^*$School of Mathematical Sciences, Queen Mary
\& Westfield College,\linebreak University of London,
Mile End Road, London E1 4NS, United Kingdom.}
\begin{abstract}
Following Isaacs (see \cite[p.~94]{isaacs}), we call a normal subgroup $N$ of a finite group $G$ \emph{large}, if $\cent_G(N) \leq N$, so that $N$ has bounded index in $G$. Our principal aim here is to establish some general results for systematically producing large subgroups in finite groups (see Theorems~A and C). We also consider the more specialised problems of finding large (non-abelian) nilpotent as well as abelian subgroups in soluble groups.
\end{abstract}
\maketitle
\section{Introduction}
Let $G$ be a finite group and suppose that $N$ is a normal subgroup of $G$ such that $\cent_G(N) \leq N$. Then we call $N$ a \emph{large} subgroup of $G$. The motivation for naming such a subgroup ``large''  stems from the observation that if $N \unlhd G$, then the factor group $G/\cent_G(N)$ is isomorphically embedded in $\Aut(N)$, and so $\left\lvert G:\cent_G(N)\right\rvert \leq \left\lvert \Aut(N)\right\rvert$. Thus, if $N$ is large,  it follows that $\left\lvert G:N\right\rvert \leq \left\lvert \Aut(N)\right\rvert$, so that $N$ has bounded index in $G$.

There are a number of examples of large subgroups in the literature. For instance, if $G$ is soluble then the (standard) Fitting subgroup $\fitt(G)$ of $G$ is large. Also, if $G$ is $\pi$-separable and $\Oh_{\pi'}(G) = 1$, then $\Oh_{\pi}(G)$ is large (this result is known as the Hall--Higman Lemma 1.2.3). Moreover, if $G$ is any finite group, then the generalised Fitting subgroup $\fast(G)$ of $G$ is large.

We briefly recall the definition of Bender's subgroup $\fast(G)$.
Call a group $G$ \emph{quasi-simple}, if $G$ is perfect and $G/\Zent(G)$ is simple.
A subnormal quasi-simple subgroup of an arbitrary group $G$ is called a \emph{component} 
of $G$ and the \emph{layer} of $G$, denoted by $\Lay(G)$, is the product of its components. 
It is known that, if $A$ and $B$ are components of $G$, then either $A = B$ or $[A, B] = 1$. 
The \emph{generalised Fitting subgroup} of $G$, denoted by $\fast(G)$,
is the product of $\Lay(G)$ and $\fitt(G)$.

We shall reserve fraktur symbols for classes of finite
groups; in particular, $\mathfrak{S}$ will denote the class of finite soluble groups, 
$\mathfrak{N}$ the class of finite nilpotent groups, and $\mathfrak{A}$ the class of finite abelian groups. 

A formation $\mathfrak{F}$ with the property that $G/\fratt(G) \in \mathfrak{F}$  for $G$ finite implies $G\in\mathfrak{F}$ is called saturated. We recall that if $\F$ is a formation, then $\F$ is called \emph{solubly saturated} provided that $G \in \F$ whenever $G/\fratt(G_{\mathfrak{S}}) \in \F$, where $G_{\mathfrak{S}}$ is the soluble radical of $G$; that is, the largest soluble normal subgroup of $G$. If $\X$ is a class of finite groups, then we shall write 
$\overline{\X}$ for its extension closure; that is, the class of all (finite) groups having a subnormal series whose quotients lie in $\X$.

\section{Main results}
Before establishing our main theorems (see Theorems~A and C below), we need to briefly discuss some auxiliary results.
\begin{lemma}\label{Lem:XBarInherit}
Let $\X$ be a class of finite groups which is closed under taking normal subgroups and quotients. Then $\overline{\X},$ the extension closure of $\X$ (i.e. the class of poly-$\X$ groups), enjoys the same closure properties.
\end{lemma}
\begin{proofof}
If $\X$ is quotient-closed then so is $\overline{\X}$ by Part (vi) of \cite[Prop. 5.4]{ResFin}. Also, a straightforward modification of the proof of Part (vii)
of \cite[Prop. 5.4]{ResFin} shows that $\overline{\X}$ inherits normal-subgroup-closure from $\X$.
\end{proofof}

\begin{lemma}\label{Lem:SimpleX}
Let $\X$ be a normal-subgroup-closed class of finite groups and let $\overline{\X}$ be its extension closure. If $G \in \overline{\X}$ and $J$ is a simple subnormal subgroup of $G,$ then $J \in \X$.
\end{lemma}
\begin{proofof}
Since $\X$ is closed under taking normal subgroups, a subnormal series for $G$ with quotients in $\X$ can be refined to yield a composition series with the same property. Since $J$ can be made the first term of a composition series for $G$, the result follows by the Jordan-H\"{o}lder theorem.
\end{proofof}

\begin{corollary}\label{Cor:NormDir}
Let $\X$ be a class of finite groups which is closed under taking normal subgroups and direct products. If $G \in \overline{\X}$ and $N$ is a minimal normal subgroup of $G,$ then $N \in \X$.
\end{corollary}
\begin{proofof}
Since $N$ is a minimal normal subgroup of $G$, we have that $N \cong J \times \cdots \times J$ for some simple group $J$. In particular, $J$ is subnormal in $G$ and thus Lemma~\ref{Lem:SimpleX} applies to yield that $J \in \X$. Now direct-product-closure of $\X$ shows that $N \in \X$.
\end{proofof}

We are now in a position to prove our first main result.

\begin{theoremA}\label{Thm:Main}
Let $\X$ be a class of finite groups which is closed under taking normal subgroups, direct products, quotients, and central extensions, and  let $\overline{\X}$ be the extension closure of $\X$. Suppose that $G \in \overline{\X},$ and that $H$ is a maximal normal $\X$-subgroup of $G$. Then $H$ is large in $G$.
\end{theoremA}
\begin{proofof}
Let $G \in \overline{\X}$, and let $H$ be a maximal normal $\X$-subgroup of $G$. We claim that $H$ is large in $G$ or, equivalently, that $\Zent(H) = \cent_G(H)$. Assume for a contradiction that $C>Z$, where $C \coloneqq \cent_G(H)$ and $Z \coloneqq \Zent(H)$.

Since $H$ is a normal subgroup of $G$, so are $C$ and $Z$, and thus $C/Z$ is a non-trivial normal subgroup of $G/Z$. Note that $G \in \overline{\X}$ implies that $G/Z \in \overline{\X}$ by Lemma~\ref{Lem:XBarInherit}. Also, $C/Z \in \overline{\X}$, again by Lemma~\ref{Lem:XBarInherit}. Now put $D/Z \coloneqq \Soc(C/Z)$. By Corollary~\ref{Cor:NormDir} we have $D/Z \in \X$ and $1<D/Z$ since $C/Z>1$ by hypothesis. Also, $H/Z \in \X$, since $\X$ is quotient-closed. Now note that $(H/Z) \cap (D/Z) = 1$ since $Z \leq H \cap D \leq H \cap C = Z$. Since $\X$ is closed under taking direct products, it follows that $HD/Z = (H/Z) \times (D/Z) \in \X$, thus central-extension-closure of $\X$ yields $HD \in \X$. However, $HD>H$, and moreover $HD$ is a normal subgroup of $G$, which contradicts the maximality of $H$.
\end{proofof}

Our next result guarantees the existence of a variety of natural classes of finite groups which are central-extension-closed.

\begin{proposition**}\label{Lem:solsat}
If $\X$ is a solubly saturated formation of finite groups that is closed under taking normal subgroups 
and such that $\mathfrak{A} \subseteq \X,$ then $\X$ is closed under taking central extensions. 
\end{proposition**}

\begin{proofof}
Let $G$ be a finite group, $1<Z \leq \Zent(G)$, and suppose that $G/Z \in \X$, with $\X$ a class of finite groups as in the statement of the proposition.

First, we argue that, since $Z$ is an abelian group, we can find an abelian group $Y$ such that $\fratt(Y) \cong Z$. This follows from the structure theorem for abelian groups, the fact that $\fratt(C_{p^{\alpha+1}}) \cong C_{p^{\alpha}}$, and \cite[Satz 6]{gasch}.

Now let $\Gamma = G \circ Y$ be the central product of $G$ with $Y$, identifying
$Z$ with the Frattini subgroup of $Y$. Then $\Gamma/Z \cong (G/Z) \times (Y/Z)$, so $\Gamma/Z \in \X$ since $Y/Z \in \A \subseteq \X$ and $\X$ is a formation thus closed under taking direct products. Also, $Z \leq \fratt(Y) \leq \fratt(\Gamma_{\mathfrak{S}})$,
so $\Gamma/\fratt(\Gamma_{\mathfrak{S}}) \in \X$ since $\X$ is quotient-closed, and thus $\Gamma \in \X$, as $\X$ is solubly saturated.  
Since $G \unlhd \Gamma$, it follows that $G \in \X$, as $\X$ is assumed to be normal-subgroup closed.
\end{proofof}

\begin{remark}
Proposition~B may also be established in a natural way by making use of Baer's local approach to 
solubly saturated formations; cf. Sec. 4 in \cite[Chap. IV]{dh}. This method may also be used to construct concrete examples of 
formations as discussed in the proposition.
\end{remark}

Our next result is a consequence of Theorem A and Proposition B.

\begin{theoremC}
Let $\X$ be a solubly saturated formation of finite groups, which is closed under taking normal subgroups, and such that $\A \subseteq \X$. Let $G\in \overline{\X}$ and let $H$ be a maximal normal $\X$-subgroup of $G$. Then $H$ is a large subgroup of $G$.
\end{theoremC}

Note that if $\X \subseteq \mathfrak{S}$, then $\X$ is solubly saturated if, and only if, it is saturated. Moreover, if $\X$ is additionally a Fitting class, then the unique maximal normal $\X$-subgroup of a finite group is its $\X$-radical.

\begin{corollary*}
\label{Cor:FittLarge}
If $G$ is a finite soluble group, then the Fitting subgroup $\fitt(G)$ of $G$ is large in $G$.
\end{corollary*}

\begin{proofof}
This follows from Theorem~C by setting $\X = \mathfrak{N}$, noting that $\overline{\mathfrak{N}} = \mathfrak{S}$, and that $\mathfrak{N}$ is a normal-subgroup-closed (solubly) saturated formation.
\end{proofof}

\begin{corollary**}
\label{Cor:GenFittLarge}
If $G$ is a finite group, then the generalised Fitting subgroup $\fast(G)$ of $G$ is large in $G$.
\end{corollary**}

\begin{proofof}
Let $\mathfrak{B}$ denote the class of finite groups which induce inner automorphisms on each of 
their chief factors. This class comprises precisely those finite groups which are \emph{quasi-nilpotent}, 
i.e., those finite groups $G$ such that $G = \fast(G)$. Since finite simple groups are quasi-nilpotent, we 
see that the extension closure of $\mathfrak{B}$ is the class of all finite groups. Further, it is known that 
$\mathfrak{B}$ is a solubly saturated Fitting formation which is closed under taking normal subgroups and that $\fast(G)$ 
is the $\mathfrak{B}$-radical of the finite group $G$; cf. \cite[p. 580]{dh}. Applying the conclusion of 
Theorem~C completes the proof.
\end{proofof}

\begin{corollary***}\label{Cor:HH123}
If $G$ is a finite $\pi$-separable group, then $\Oh_{\pi', \pi}(G)$ is large in $G$.
\end{corollary***}

\begin{proofof}
Consider the class of groups which possess a normal Hall $\pi'$-subgroup. This class
is easily seen to be a saturated Fitting formation which contains all abelian groups, and whose 
extension closure is the class of $\pi$-separable groups. Note that if $G$ is $\pi$-separable
then $\Oh_{\pi', \pi}(G)$ is the radical associated with the class mentioned above. Thus
Theorem~C applies and gives what we want.
\end{proofof}

It is easy to see that Corollary~F is just an equivalent formulation of \cite[Lemma 1.2.3]{psoluble};
cf. also \cite[Thm. 3.21]{isaacs}.

\section{Large nilpotent subgroups in soluble groups}
Our next result, which generalises an observation communicated to us by Marty Isaacs, concerns existence of large subgroups in finite soluble groups, which are of bounded nilpotency class. We note that Proposition~G does not follow from our previous results, since the class  $\X = \mathfrak{N}_c$ of finite nilpotent groups of class at most $c$ relevant here is not closed under central extensions. 

\begin{proposition*****}
Suppose that $G$ is finite and soluble, and let $H \unlhd G$ be maximal in $G$ 
with the property that it is nilpotent of class at most $c,$ where $c \geq 2$ is a given integer. Then $H$ is large in $G$.
\end{proposition*****}
\begin{proofof}
Assume that the assertion is false, and let $(G,H)$ be a counterexample.
Write $C \coloneqq \cent_G(H)$ and $Z \coloneqq \Zent(H)$. 
Let $1 = C_0 < \cdots < C_r = C$ be the derived series of $C$. Since $C$ is not 
contained in $H$, there exists a smallest positive index $j$ such that 
$C_j \nleqslant H$. In this situation, we claim that 
\begin{equation}\label{Eq:MainIneq}
L_i(HC_j) \leq L_i(H), \quad \text{for} \,\ i \geq 2,
\end{equation}
where $(L_i(K))_{i\geq 0}$ denotes the lower central series of the group $K$. We first note that 
\begin{equation}
L_1(HC_j) = (HC_j)' = H'C_{j-1} \leq H'(C \cap H) = H'\Zent(H),
\end{equation}
where we have used the minimality of the index $j$ to deduce that $C_{j-1} \leq H$. Thus, 
using the fact that $C$ and $Z$ are normal in $G$ (because $H$ is), we get
\begin{multline*}
L_2(HC_j) = [H C_j, L_1(H C_j)] \leq [HC_j, H' \Zent(H)] \\
= [H, H'] \cdot [H, \Zent(H)] \cdot [C_j, H'] \cdot [C_j, \Zent(H)] 
= [H, H'] = L_2(H). 
\end{multline*}
Assuming inductively that 
$L_i(H C_j) \leq L_i(H)$ for some $i\geq 2$, we now find that
\begin{multline*}
L_{i+1}(H C_j) = [H C_j, L_i(H C_j)] \leq [H C_j, L_i(H)] \\
= [H, L_i(H)] \cdot [C_j, L_i(H)] = [H, L_i(H)] = L_{i+1}(H), 
\end{multline*}
whence \eqref{Eq:MainIneq}. 

Next, we claim that 
\begin{equation}
\label{Eq:MainIneq2}
c(HC_j) \leq c(H) + 1,
\end{equation}
where $c(K)$ denotes the nilpotency class of $K$.  In order to see this, suppose first that $c(H) \leq  1$, i.e., that $H$ is abelian. Then we need to show that $c(H C_j) \leq 2$ or, equivalently, that $(H C_j)^\prime \leq \Zent(H C_j)$. However, we have 
\[
(H C_j)^\prime = H^\prime C_{j-1} = C_{j-1} \leq H\cap C = \Zent(H) \leq \Zent(H C_j),
\]
so that our claim \eqref{Eq:MainIneq2} holds in this case. Next, let $c(H) \geq 2$. Then we may use (\ref{Eq:MainIneq}) to deduce that 
\[
c(H C_j) \leq c(H) \leq c(H) + 1,
\]
as required, whence \eqref{Eq:MainIneq2}. 

We are now in a position to finish the proof. Since $C_j \not\leq H$, we have $H C_j > H$, and also $H C_j \unlhd G$, since $H \unlhd G$ and $C_j$ is characteristic in $C$, which in turn is normal in $G$. Moreover, in view of \eqref{Eq:MainIneq2}, if $c(H)< c$, we have $c(H C_j) \leq c$, contradicting the maximality of $H$. Hence, we may assume that $c(H) = c \geq 2$. In this situation however, we may apply \eqref{Eq:MainIneq} to conclude again that $c(H C_j) \leq c$.  This final contradiction to the maximality of $H$ shows that no counterexample to our claim exists, and the proof is complete. 
\end{proofof} 
\begin{remark}
A result analogous to Proposition~G holds for the class $\mathfrak{S}_d$ of finite soluble groups of derived length at most $d$, where $d \geq 2$. The proof is similar to that of Proposition G, but is considerably easier due to the fact that only commutators play a role.
\end{remark}

\section{Large abelian subgroups in soluble groups}
We note that there is no analogue of Proposition G for maximal 
\emph{abelian} normal subgroups of soluble groups. Of course, \emph{some} soluble 
groups (like $A_4$) contain maximal abelian normal subgroups which are large 
(i.e. self-centralising), but not all of them do. In $\mathrm{SL}_2(3)$, the centre is a 
maximal abelian normal subgroup, but, of course, it does not contain its centraliser, 
which is the whole group. Our final result exhibits a class of finite soluble groups properly containing the class of 
finite supersoluble groups, for which an analogue 
of Proposition G does exist. 

Denote by $\X_0$ the class of those finite soluble groups $G$ whose supersoluble residual is either trivial 
or a minimal normal subgroup of $G$. Note that $A_4\in \X_0$ but that $A_4$ is not supersoluble; in particular, $\X_0$ is strictly larger than the class of all finite supersoluble groups. 

\begin{proposition******}
If $G\in\X_0,$ then every maximal abelian normal subgroup of $G$ is large \nolinebreak in \nolinebreak $G$.
\end{proposition******}
\begin{proofof}
First, we treat the case where $G$ is supersoluble. Let $H$ be a maximal abelian normal subgroup of 
the supersoluble group $G$, and write $C \coloneqq \cent_G(H)$. As $H \unlhd G$, so $C \unlhd G$ 
and since $H$ is abelian we have $H \leq C$. We claim that, in fact, $H=C$. We assume that this assertion 
is false, so that $H<C$, and work to derive a contradiction.

Consider the non-trivial group $C/H$, which is normal in the supersoluble group $G/H$. Since every non-trivial normal subgroup of a supersoluble group contains a non-trivial cyclic normal subgroup, it follows that there is an element $x \in C\setminus H$ such that $\langle xH \rangle = \langle 
H,x \rangle/H \unlhd G/H$.

Now, consider the group $\langle H,x \rangle$. Since $x$ centralises $H$ and $H$ is abelian, it follows that $\langle H,x \rangle$ is abelian, and it is normal in $G$ since it is the full preimage of a normal subgroup of $G/H$. Moreover, $\langle H,x \rangle > H$ since $x \in C\setminus H$, contradicting the fact that $H$ was chosen maximal with respect to being abelian and normal in $G$.

Next, let $G \in \X_0$ be such that its supersoluble residual $S$ is a minimal normal subgroup of $G$; in particular, $S$ is abelian since $G$ is soluble. Let $A$ be a maximal abelian normal subgroup of $G$. Then we have $S \leq A$, since otherwise,
by minimality of $S$, $A \cap S =1$ and thus $AS = A \times S$ is abelian (and normal in $G$), contradicting
the maximality of $A$.

Now we argue as before. Suppose, for a contradiction, that $C \coloneqq \cent_G(A) >A$. Since $A \geq S$, it follows that $G/A$ is
isomorphic to a quotient of the supersoluble group $G/S$, thus is itself supersoluble. Then $C/A$ is a non-trivial normal subgroup
of $G/A$, hence contains some non-trivial normal cyclic subgroup, say $K/A$. As $A$ is central in $C$, so $A$ is central in $K$, and thus
$K$ is normal in $G$ and abelian, since it is central-by-cyclic. This again contradicts the maximality of $A$, thus we have $C=A$, as desired.
\end{proofof}

The property of the class $\X_0$ expressed in Proposition~$H$ does not characterise $\X_0$, however, in that there exist soluble groups having every maximal abelian normal subgroup self-centralising (i.e. large), which are not contained in $\X_0$.  Indeed, as GAP \cite{GAP4} informs, the two smallest groups $G$ such that $G\not \in \X_0$, while every maximal abelian normal subgroup is large in $G$, are of order $48$ (SmallGroup[48, 3] and SmallGroup[48,50]). A larger example of this sort is the group $G= A_4\times A_4$. Here, the supersoluble residual is the $2$-Sylow subgroup of $G$, while the only minimal normal subgroups of $G$ are the $2$-Sylow subgroups of the individual factors. It would be of some interest to describe those groups in $\mathfrak{S} \setminus \X_0$ for which every maximal abelian normal subgroup is large.

\section*{Acknowledgments}
The authors would like to express their gratitude to Marty Isaacs for some stimulating discussions on the topic of the present paper.

The first author acknowledges financial support in the form of a fellowship from the Lady Davis Foundation
at the Hebrew University of Jerusalem, and grants from the European Research Council and the Israel Science 
Foundation no. 1117/13.

\end{document}